\documentclass[11 pt]{amsart}  

\usepackage{graphicx}
\usepackage{amsmath,amsthm,amscd,amssymb}
\usepackage{amsfonts,color}
\usepackage{latexsym}
\usepackage{epsfig,epstopdf}
\usepackage{enumerate}

\newcommand{\ds}{\displaystyle}

\newcommand{\R}{\mathbb{R}}            

\newcommand{\Z}{\mathcal{Z}}

\newcommand{\U}{\mathcal{U}}

\newcommand\norm[1]{\left\| #1\right\|}
\newcommand\inner[1]{\langle #1\rangle}

\newcommand{\Aa}{\mathbb{A}}
\newcommand{\Ba}{\mathbb{B}}
\newcommand{\Ia}{\mathbb{I}}
\newcommand{\Ga}{\mathbb{M}}
\newcommand{\Fa}{\mathbb{F}}

\newtheorem{teo}{\sc Theorem}[section]

\newtheorem{lemma}{\sc Lemma}[section]

\newcommand{\demo}[1]{\noindent {\bf Proof of Theorem #1.} }
\newcommand{\findem}{\mbox{\framebox [1ex]{}}}
\newcommand{\dem}{\noindent {\bf Proof} \mbox{}}    

\begin{document}

\title[Controllability of  Semilinear Equations]
{Controllability of Impulsive Semilinear Evolution Equations with Memory and Delay in Hilbert Spaces}
\date{\today}
\author[C. GUEVARA  AND H. LEIVA ]{C. Guevara$^1$ and H. Leiva$^2$  }
\address{$^{1}$ Triptych,LLC\\Gilbert AZ 85296\\ AND \\Arizona State University \\
        Simon A. Levin Mathematical, Computational and Modeling Sciences Center \\
          Tempe, AZ 85281-USA} \email{cristi.guevara@asu.edu, Triptych.LLC.AZ@gmail.com}
   \address{$^{2}$ School of Mathematical Sciences and Information Technology \\
           Universidad Yachay Tech \\
         San Miguel de Urcuqui, Ecuador}
          \email{hleiva@yachaytech.edu.ec, hleiva@ula.ve}

\thanks{$\dagger$ This work has been supported by YachayTech University }

\begin{abstract}
Inspired in our work on the controllability for the semilinear with memory   \cite{Carrasco-Guevara-Leiva:2017aa,  Guevara-Leiva:2016aa, Guevara-Leiva:2017aa}, we present the general cases for the approximate controllability of impulsive semilinear evolution equations in a Hilbert space with memory and delay terms which arise from reaction-diffusion models. We prove that, for each initial and an arbitrary neighborhood of a final state, one can steer the system from the initial condition to this neighborhood of the final condition with an appropriated  collection of admissible controls thanks to the  delays. Our proof is based on semigroup theory and  A.E. Bashirov et al. technique \cite{
Bashirov-Ghahramanlou:2015aa, Bashirov-Jneid:2013aa, Bashirov-Mahmudov:2007aa} which avoids fixed point theorems.
\end{abstract}

\maketitle


\section{Introduction}

Control theory theory has gained a huge scope in the modern times. It comes very handy to understand, diagnose, modify, improve and predict behavior of biological processes, engineering practices, robotics, internet of things (IoTs),  and dynamics of business, social, and political systems among others.

Today, in our data-driven society control theory is of vital importance, for instance to monitor the stock market; every entity to be controlled can be associated to a controller and every data that is being monitored becomes the data point, and thus control laws applied.  In general,   control theory addresses how a systems can be modified through feedback, in particular, how an arbitrary initial state can be directed either exactly or approximated close to a given final state using a set of admissible controls. Furthermore, the controllability is a robust property, in the sense, if a system is controllable, all modes of the system can be perturbed from the input and inherent phenomena such as abrupt changes, delays and dependance on prior behavior would not modify the controllability of the system.  Thus, the conjecture is that controllability of a system will not change due to perturbations corresponding to delays, impulses or some type of memories.

For the purpose of this paper, consider the following class of impulsive semilinear evolution equations with memory and delay on Hilbert spaces $\U$ and $\Z$
\begin{equation}\label{eq:class}
\left\{
\begin{array}{lr}
z' =  -\Aa z+ \Ba u + a\; \ds \int_{0}^{t}\Ga(t,s ,z_{s})ds + b \; \Fa(t,z_{t},u(s)) ,& z\in Z^{1},\; t\geq 0, \\
z(s)  = \Phi(s), & s \in [-r,0],\\
z(t_{k}^{+})  = z(t_{k}^{-})+\Ia_{k}(t_k, z(t_{k}),u(t_{k})), &  k=1,2,3, \dots, p,
\end{array}
\right.
\end{equation}
such that  $t_k \in (0, \tau)$ with $t_k<t_{k+1}$,  $a,b \in \R$, the operator $\Aa :D(\Aa) \subset \Z \rightarrow \Z$ is sectorial and $-\Aa$ is the infinitesimal generator of a compact analytic semigroup of uniformly bounded linear operators $\{T(t) \}_{t \geq 0} \subset \Z$, with $0 \in \rho(\Aa)$. Therefore, fractional power operators $\Aa^{\beta}, \quad 0 < \beta \leq 1$, are well defined. And since $A^{\beta}$ is a closed operator, its domain $D(A^{\beta})$ is a Banach space endowed with the graph norm
$$
\|z\|_{\beta} = \|\Aa^{\beta}z\|, \ \ z \in D(\Aa^{\beta}).
$$
This Banach space is denoted by $\Z^{\beta}=D(\Aa^{\beta})$ and it is dense in $\Z$. The standard notation $z_{t}$ defines a function from $[-r,0]$ to $\Z^{\beta}$ by $z_{t}(s) = z(t+s), -r \leq  s \leq 0$ for fixed $0 < \beta \leq 1$, and initial state $\Phi \in \mathcal{PW}= \mathcal{PW}\left(-r,0;  \Z^{\beta}\right)$, with $\mathcal{PW}$ the space of piecewise continuous functions.
Here,   $u \in L^{2}(0,\tau; \U)$ represents the control, $\Ba: \U \longrightarrow \Z$ is a bounded linear operator  and $\Ia_{k}:[0, \tau]\times \Z^{\beta} \times \U \to \Z,$\; $\Fa: [0, \tau] \times \mathcal{PW}(-r,0;  \Z^{\beta}) \times \U \to \Z$  and $\Ga: [0,\tau]^2\times \mathcal{PW}(-r,0;  \Z^{\beta}) \to \Z$ are smooth enough functions.\\
Moreover, in order to prove the controllability of the corresponding linear system, we shall assume that the strongly continuous semigroup $\{T(t) \}_{t \geq 0}$ generated by $\Aa$ satisfies the following spectral decomposition:
$$
 \Aa z = \sum_{j = 1}^{\infty} \lambda_{j} \sum_{k = 1}^{\gamma_j} <z, \phi_{j,k}> \phi_{j,k}, \ \ z \in \Z
$$
with the eigenvalues $0 < \lambda_1 < \lambda_2 < \cdots < \cdots \lambda_n \to \infty$ of $\Aa$ having finite multiplicity $\gamma_j$
equal to the dimension of the corresponding eigenspaces, and $\{\phi_{j,k}\}$ is a complete orthonormal set of
eigenfunctions of $A$. So, the  strongly continuous semigroup $\{T(t) \}_{t \geq 0}$ generated by  $-\Aa$ is given by
$$
T(t)z =  \sum_{j = 1}^{\infty} e^{-\lambda_{j}t}\sum_{k =1}^{\gamma_j} <z, \phi_{j,k}> \phi_{j,k}, \ \ z \in \Z,
$$
and
$$
 \Aa^{\beta} z = \sum_{j = 1}^{\infty} \lambda_{j}^{\beta} \sum_{k = 1}^{\gamma_j} <z, \phi_{j,k}> \phi_{j,k}, \ \ z \in \Z^{\beta}
$$
As a consequence, we have the following estimate:
$$
\parallel T(t)\parallel\leq Ke^{-\mu t},\quad t\geq0, \ \ \mu >0.
$$
Also, we will assume the following hypothesis\\
(H1) $\Ba^{*}\phi_{j,k}$ are linearly independent in $\Z$.


Assuming that $\Ga \in L^{\infty}( [0,\tau]^2 \times \mathcal{PW}(-r,0;  \Z^{\beta})$ and  $ \Fa, \Ia_{k}$  are smooth enough so that for every $\Phi \in \mathcal{PW}([-r,0],\Z^{\beta})$ and every control $u$, the equation \eqref{eq:class} admits only one mild solution on  $0\leq t_{0}\leq t\leq \tau$ given by
\begin{align}\label{eq:mild-class}
z(t) =\; & \ds T(t)\Phi(0)+\int_{0}^{t}T(t-s)\left[\Ba u(s)+a \left(\int_{0}^{s}\Ga(s,l,z_{l})dl\right)\right]ds \\
& \;+   b \ds \int_{0}^{t}T(t-s)\Fa(s,z_{s},u(s))ds  +   \sum_{0 < t_k < t} T(t-t_k )\Ia_{k}(t_k,z(t_k), u(t_k)).  \nonumber
  \end{align}
Additionally, suppose  there exist a continuous function $\rho : \mathbb{R}_{+} \rightarrow [0, \infty )$ such that, for all $(t, \Phi,u) \in [0, \tau] \times \mathcal{PW}(-r,0;  \Z^{\beta}) \times \U$ the following inequality holds
\begin{equation}\label{eq:bound}
\norm{\Fa(t,\Phi,u)}_{\Z}  \leq    \rho(\| \Phi(-r) \|_{\Z}).
\end{equation}
In particular, $\rho(\xi)=e(\xi)^\alpha + \eta$, with $\alpha \geq 1$.

Recall that \eqref{eq:mild-class}  is said to be {\bf approximate controllable}  on $[t_{0},\tau]$ if for every $z_0$,
$z_1\in\Z$, and $\epsilon>0$ there exists $u\in L^{2}(t_{0},\tau;\U)$ such
that the mild solution $z(t)$ corresponding to $u$
verifies:
\begin{equation}\label{eq:goal}
\norm{ z(\tau) - z_{1}}_{\Z}<\epsilon.
\end{equation}
The purpose of this paper is to prove
\begin{teo}\label{teo:main}
Under the above assumptions the  semilinear evolution equation with memory, delay and impulses  \eqref{eq:class} is approximately controllable on $[0, \tau]$.
\end{teo}


\section{Linear System characterization}\label{sec:lineal}
 In order to prove  Theorem \ref{teo:main}, in this section we shall prove the approximate controllability of the following linear evolution equation without memory, impulses and delay
\begin{equation}\label{eq:linear}
\left\{%
\begin{array}{lll}
    z'(t) = -\mathbb{A} z(t) + \Ba u(t),\\
    z(t_{0}) = z_{0},
    \end{array}%
\right.
\end{equation}
where $z_{0}\in \Z$ and $u\in L^{2}([0,\tau];\U)$, obtained by setting $a=b=0$ in \eqref{eq:class}.\\
For the system \eqref{eq:linear} and $\tau>0$, we have the following notions:
\begin{enumerate}
\item $G_{\tau\delta}$ is the controllability  operator defined by
\begin{eqnarray*}
G_{\tau\delta}: L^2(\tau-\delta,\tau;\U) \longrightarrow& \Z\\
u\longmapsto&\ds \int_{\tau-\delta}^{\tau}T(\tau-s)\Ba u(s)ds,
\end{eqnarray*}
with corresponding adjoint $G^*_{\tau\delta}$ given by
\begin{eqnarray*}
G^*_{\tau\delta}:  \Z \longrightarrow& L^2(\tau-\delta,\tau;\U)\\
z\longmapsto& \Ba ^{*}T^{*}(\tau-\cdot)z.
\end{eqnarray*}
\item The Gramian controllability operator is
\begin{equation*}
Q_{\tau \delta*} = G_{\tau\delta}G_{\tau\delta}^{*}= \int_{\tau-\delta}^{\tau}T(\tau-t)\Ba \Ba ^{*}T^{*}(\tau-t)dt.
\end{equation*}
\end{enumerate}

In general, for linear bounded operator $G$ between Hilbert spaces $\mathcal{W}$ and $\mathcal{Z}$, the following lemma holds (see \cite{Bashirov-Kerimov:1997aa,Bashirov-Mahmudov:1999aa, Leiva-Merentes:2013aa}).
\begin{lemma}\label{Lema1}
The approximate controllability of \eqref{eq:linear}  on $[\tau-\delta,\tau]$ is equivalent to  any of the following statements
\begin{enumerate}[(a)]
\item $\overline{Rang(G_{\tau\delta})}=\Z.$
\item $\Ba ^{*}T^{*}(\tau-t)z={0}, \quad t \in [\tau-\delta,\tau]  \Rightarrow z =0.$
\item For $0\neq z \in\ Z^{1},  \  \ \inner{ Q_{\tau\delta}z,z}>0.$
\end{enumerate}
\end{lemma}
Theorem \ref{A1.5}   and Lemma \ref{lema} characterized the controllability of the system \eqref{eq:linear}, their proofs and details can be found in \cite{Bashirov-Kerimov:1997aa,Bashirov-Mahmudov:1999aa, Curtain-Pritchard:1978aa, Curtain-Zwart:1995aa, Leiva-Merentes:2013aa}\\
Under the hypothesis (H1), we can prove the approximate controllability of linear system \eqref{eq:eq:linear}
\begin{teo}\label{main1}
If  vectors $\Ba^{*}\phi_{j,k}$ are linearly independent in $\Z$, then
the system \eqref{eq:eq:linear} is approximately controllable on $[\tau-\delta,\tau]$, for $0 < \delta < \tau$.
\end{teo}
\demo We shall apply condition (b) from the foregoing Lemma. In fact, clearly that $T^{*}(t) = T(t)$, and suppose that:
$$
\Ba ^{*}T^{*}(\tau-t)z={0}, \quad t \in [\tau-\delta,\tau].
$$
i.e.,
$$
\sum_{j = 1}^{\infty} e^{-\lambda_{j}(\tau-t)}\sum_{k =1}^{\gamma_j} <z, \phi_{j,k}> \Ba ^{*} \phi_{j,k}=0, \quad t \in [\tau-\delta,\tau].
$$
i.e.,
$$
\sum_{j = 1}^{\infty} e^{-\lambda_{j}t}\Ba ^{*}\sum_{k =1}^{\gamma_j} <z, \phi_{j,k}> \phi_{j,k}=0, \quad t \in [0, \delta].
$$
From Lemma 3.14 from \cite{Curtain-Pritchard:1978aa}, we get that
$$
\sum_{k =1}^{\gamma_j} <z, \phi_{j,k}> \Ba ^{*} \phi_{j,k}=0
$$
Now, from the hypothesis (H1), we get that $<z, \phi_{j,k}> =0, \quad k= 1,2, \dots, \gamma_j ; j=1,2,3, \dots.$ \\ Since $\{\phi_{j,k}\}$ is a complete orthonormal set on $\Z$, we conclude that $z=0$. This completes the proof of the approximate controllability of the linear system \eqref{eq:linear}.

\hfill \findem

Another characterization of the approximate controllability of system \eqref{eq:linear} follows from Lemma \ref{Lema1}:

 \begin{teo}\label{A1.5} The system \eqref{eq:linear} is approximately
controllable on $[\tau-\delta,\tau]$,  for $0 < \delta < \tau$, if and only if any one of the following conditions hold:
\begin{enumerate}
\item $\ds\lim_{\alpha \to 0^+} \alpha(\alpha I +Q_{\tau\delta}^{*})^{-1}z =0 $.\\
\item If $z\in Z^{1}$,  $0<\alpha \leq 1$ and $u_{\alpha}=G_{\tau\delta}^{*}(\alpha I +
Q_{\tau\delta}^{*})^{-1}z$, then
$$
G_{\tau\delta}u_{\alpha}=z - \alpha(\alpha I+ Q_{\tau\delta})^{-1}z \quad \mbox{and} \quad
\displaystyle\lim_{\alpha\to 0}G_{\tau\delta}u_{\alpha}=z.
$$
Moreover, for each $v\in L^{2}([\tau-\delta,\tau];\U)$, the sequence of controls
$$
u_{\alpha}=G_{\tau\delta}^{*}(\alpha I +
Q_{\tau\delta}^{*})^{-1}z + (v-G_{\tau\delta}^{*}(\alpha I +
Q_{\tau\delta}^{*})^{-1}G_{\tau\delta}v),
$$
satisfies
$$
G_{\tau\delta}u_{\alpha}=z-\alpha(\alpha I +
Q_{\tau\delta}^{*})^{-1}(z-G_{\tau\delta}v)
\quad
\mbox{and}
\quad
\displaystyle\lim_{\alpha\to 0}G_{\tau\delta}u_{\alpha}=z,
$$
with the error $E_{\tau\delta}z=\alpha(\alpha I +
Q_{\tau\delta})^{-1}(z+G_{\tau\delta}v),\;\alpha\in(0,1].
$
\end{enumerate}
\end{teo}
Theorem \ref{A1.5} indicates that the family of linear operators $
\Gamma_{\tau\delta}=G_{\tau\delta}^{*}(\alpha I +
Q_{\tau\delta}^{*})^{-1}
$
satisfies the following limit
$$
\displaystyle\lim_{\alpha\longrightarrow 0}G_{\tau\delta}\Gamma_{\tau\delta}=I,
$$
in the strong topology.

\begin{lemma}\label{lema}
$Q_{\tau\delta}> 0$, if and only if, the linear system \eqref{eq:linear} is controllable on $[\tau-\delta, \tau]$.
Moreover, for  given initial state $y_0$ and  final state $z_{1}$, there exists a sequence of controls $\{u_{\alpha}^{\delta}\}_{0 <\alpha \leq 1}$ in the space $L^2(\tau-\delta,\tau;\U)$, defined by
$$
u_{\alpha}=u_{\alpha}^{\delta}= G_{\tau\delta}^{*}(\alpha I+ G_{\tau\delta}G_{\tau\delta}^{*})^{-1}(z_{1} - T(\tau)y_0),
$$
such that the solutions $y(t)=y(t,\tau-\delta, y_0, u_{\alpha}^{\delta})$ of the initial value problem
\begin{equation}\label{IVL}
\left\{
\begin{array}{l}
y'=\Aa y+\Ba u_{\alpha}(t), \ \  y \in \Z, \ \ t>0,\\
y(\tau-\delta) = y_0,
\end{array}
\right.
\end{equation}
satisfies
\begin{equation}\label{eq:limit}
\lim_{\alpha \to 0^{+}}y(\tau)
= \lim_{\alpha \to 0^{+}}\left(T(\delta)y_0 + \int_{\tau-\delta}^{\tau}T(\tau-s)\Ba u_{\alpha}(s)ds \right)= z_{1}.
\end{equation}
\end{lemma}

\section{Controllability of the Semilinear System}\label{sec:semilineal}

This section is devoted to prove the main result of this paper, the approximate controllability of the impulsive semilinear evolution equations in a Hilbert space with memory and delay terms \eqref{eq:class}: Theorem \ref{teo:main}.
The approach to obtain \eqref{eq:goal} consist in construct a sequence of controls conducting the system from the initial condition $\Phi$ to a small ball around $z_1,$ taking advantage of the delay, which  allows us to pullback the corresponding family of  solutions to a fixed trajectory in short time interval. Now, we are ready to present the proof of our main result

\dem{\ref{teo:main}} Let $\epsilon>0$, and given $\Phi\in \mathcal{PW}$ and a final state $z_{1} \in \Z$.  Consider any $u\in L^{2}([0,\tau];\U)$ and the corresponding mild solution \eqref{eq:mild-class}  of the initial value problem \eqref{eq:class}, denoted by $z(t)=z(t,0,\Phi,u)$. For $0\leq\alpha \leq 1,$ define the control $u_{\alpha}^{\delta}\in L^{2}([0,\tau];\U)$  as follows
$$
u_{\alpha}^{\delta}(t)=\left\{\begin{array}{ccl}
                         u(t), &&0\leq t\leq \tau-\delta, \\
                         u_{\alpha}(t), &\quad& \tau-\delta\leq t\leq \tau,
                       \end{array}\right.
$$
with $
u_{\alpha}= \Ba^{*}T^{*}(\tau-t)(\alpha I+ G_{\tau\delta}G_{\tau\delta}^{*})^{-1}(z_{1} - T(\delta)z(\tau-\delta)).
$
Thus,

 $0<\delta<\tau-t_{p}$ and the corresponding mild solution at time $\tau$ can be written as follows:
 \begin{eqnarray*}
 \ds
z^{\delta,\alpha}(\tau) &=&
 \ds T(\tau)\Phi(0) +\int_{0}^{\tau}T(\tau-s)
	\left[ \Ba   u_{\alpha}^{\delta} (s)
		+ a \int_0^s \Ga(s,l,z^{\delta,\alpha}_{l})dl\right]ds+ \\
	&&+b  \int_{0}^{\tau}T(\tau-s)\Fa(s,z^{\delta,\alpha}_{s},u_{\alpha}^{\delta}(s))ds+
\sum_{0 < t_k < \tau} T(t-t_k )\Ia_{k}(t_k,z^{\delta,\alpha}(t_k), u_{\alpha}^{\delta}(t_k))\\
&=&T(\delta)\left\{T(\tau-\delta)\Phi(0)
+\int_{0}^{\tau-\delta}T(\tau-\delta-s) \left(\Ba   u_{\alpha}^{\delta} (s)+b\; \Fa(s,z^{\delta,\alpha}_{s},u_{\alpha}^{\delta}(s))\right)ds\right.\\
&&\qquad\quad+a\int_{0}^{\tau-\delta}T(\tau-\delta-s)  \int_0^s \Ga(s,l, z^{\delta,\alpha}_{l})dlds\\
&&\qquad\quad\left.+ \sum_{0 < t_k < \tau-\delta} T(t-\delta-t_k )\Ia_{k}(t_k,z^{\delta,\alpha}(t_k), u_{\alpha}^{\delta}(t_k))\right\}+\\
&& + \int_{\tau-\delta}^{\tau}T(\tau-s)\left(\Ba u_{\alpha}(s)+
b\;\Fa(s,z^{\delta,\alpha}_{s},u_{\alpha}^{\delta}(s))+a\;\int_0^s\Ga(s,l,z^{\delta,\alpha}_{l})dl\right)ds.
\end{eqnarray*}
Thus,
\begin{align*}
z^{\delta,\alpha}(\tau)  = &T(\delta)z(\tau-\delta)+ \int_{\tau-\delta}^{\tau}T(\tau-s)\left(\Ba u_{\alpha}(s)+ b\;
\Fa(s,z^{\delta,\alpha}_{s},u_{\alpha}^{\delta}(s))\right)ds \\
&+a\; \int_{\tau-\delta}^{\tau}T(\tau-s)\int_0^s\Ga(s,l,z^{\delta,\alpha}_{l})dlds.
\end{align*}
Note that the corresponding solution $y^{\delta,\alpha}(t)=y(t,\tau-\delta,z(\tau-\delta),u_{\alpha})$ of the initial value problem \eqref{IVL} at time $\tau$ is:
$$
y^{\delta,\alpha}(\tau)=T(\delta)z(\tau-\delta)+ \int_{\tau-\delta}^{\tau}T(\tau-s)\Ba_{\varpi} u_{\alpha}(s)d.,
$$
Hence,
$$
z^{\delta,\alpha}(\tau)-y^{\delta,\alpha}(\tau)=
\int_{\tau-\delta}^{\tau}T(\tau-s)\left(\int_{0}^{s}b\; \Fa(s,z^{\delta,\alpha}_{s},u_{\alpha}^{\delta}(s)+a\;\Ga(s,l,z^{\delta,\alpha}_{l})dl)\right)ds,
$$
 and from condition (\ref{eq:bound})we obtain that
\begin{align*}
  \norm{ z^{\delta,\alpha}(\tau)-y^{\delta,\alpha}(\tau)} & \leq |b|\int_{\tau-\delta}^{\tau} \norm{ T(\tau-s)}\rho\left( \norm{z^{\delta,\alpha}(s-r)}\;\right)ds \\
   & + |a| \int_{\tau-\delta}^{\tau}\norm{ T(\tau-s)}\int_{0}^{s} \norm{\Ga(s,l,z^{\delta,\alpha}_{l})}dlds.
\end{align*}
Observe that
 $0< \delta< r$ and $\tau-\delta \leq s\leq \tau$, thus $$l-r \leq s-r \leq \tau-r< \tau-\delta.
 $$
Therefore,
$
z^{\delta,\alpha}(l-r)=z(l-r) $  and $ z^{\delta,\alpha}(s-r)=z(s-r),
$ implying that for $\epsilon>0$ there exists $\delta>0$ such that
\begin{align*}
 \norm{z^{\delta,\alpha}(\tau)-y^{\delta,\alpha}(\tau)} & \leq|b| \int_{\tau-\delta}^{\tau}\norm{ T(\tau-s)}\rho\left(\norm{z(s-r)}\;\right)ds \\
  &\quad +|a|  \int_{\tau-\delta}^{\tau}\norm{T(\tau-s)}\int_{0}^{s}\norm{ \Ga(s,l,z_{l})} dlds  \\
   & <  \displaystyle\frac{\epsilon}{2}.
\end{align*}
Additionally, for $0<\alpha <1$, Lemma \ref{lema}  \eqref{eq:limit} yields
$$
 \norm{ y^{\delta,\alpha}(\tau)-z_{1}}  <  \frac{\epsilon}{2}.
$$
Thus,
$$
\begin{array}{lll}
 \norm{ z^{\delta,\alpha}(\tau)-z_{1}}  & \leq &  \norm{ z^{\delta,\alpha}(\tau)-y^{\delta,\alpha}(\tau)} +  \norm{ y^{\delta,\alpha}(\tau)-z_{1}}  <  \frac{\epsilon}{2}+ \frac{\epsilon}{2}=\epsilon,
\end{array}
$$
which completes our proof.

\hfill \findem

\section{Final Remarks}\label{final}
J.P. LaSalle in \cite{LaSalle:1960aa} wrote the following ``It is never possible to start the system exactly in its equilibrium state, and
the system is always subject to outside forces not taken into account by the differential equations. The system is disturbed and is displaced slightly from its equilibrium state. What happens? Does it remain near the equilibrium state? This is stability. Does it remain near the equilibrium state and in addition tend to return to the equilibrium? This is asymptotic stability". Continuing with what J.P. LaSalle said, we conjecture that real life systems are always under the influence of impulses, delays, memory, nonlocal conditions and noises, which are intrinsic phenomena no taken into account  by the mathematical model that is represented by a differential equation, and if we consider these intrinsic phenomena as model disturbances, we have noticed that controllability is not lost, in other words, under certain conditions, these phenomena, seen as system disturbances do not destroy its controllability. So, we can add impulses, delays, non-local conditions, noise, etc. and some properties of the system persist.

%

\end{document}